\documentclass[a4paper,10pt,fleqn]{article}
\usepackage{epsfig,epsf,graphics,graphicx,german,amssymb,amsmath,latexsym,color,psfrag}
\usepackage{hyperref}
\usepackage[authoryear]{natbib}
\usepackage[english]{babel}
\usepackage{lscape}
\newcommand{\gl}{\;\:=\;\:}

\DeclareMathOperator{\MutInf}{MutInf}
\DeclareMathOperator{\sMutInf}{sMutInf}
\DeclareMathOperator{\sign}{sign}
\DeclareMathOperator{\HS}{HS}

\newcommand{\ve}{\varepsilon}

\begin{document}
\begin{center}
{\bf\huge Comparing measures of association in 2x2 probability tables} \\
\vspace{3cm}
{\bf Dirk Hasenclever$^1$ and Markus Scholz$^{1,2}$} \\
$^1$Institute for Medical Informatics, Statistics and Epidemiology \\
$^2$LIFE Research Center \\
University of Leipzig 
\end{center}

\vspace{1cm}\noindent
{\bf Running head:} Comparing measures of association \\
\\ \\ \\ \noindent
{\bf Corresponding author:} 
Dirk Hasenclever \\
University of Leipzig \\
Institute for Medical Informatics, Statistics and Epidemiology \\
Haertelstrasse 16-18\\
04107 Leipzig \\
Germany \\
Telephone: +49 341 97 16121 \\
Fax: +49 341 97 16109 \\
E-Mail: dirk.hasenclever@imise.uni-leipzig.de

\newpage
\section*{Abstract}
Measures of association play a role in selecting 2x2 tables exhibiting strong dependence in high-dimensional binary data. 
Several measures are in use differing on specific tables and in their dependence on the margins. 

We study a 2-dimensional group of margin transformations on the 3-dimensional manifold T of all 2x2 probability tables. The margin transformations allow introducing natural coordinates 
that identify T with the real 3-space such that the x-axis corresponds to log(sqrt(odds-ratio)) and margins vary on planes x=const. We use these coordinates to visualise and compare 
measures of association with respect to their dependence on the margins given the odds-ratio, their limit behaviour when cells approach zero and their weighting properties. 

We propose a novel measure of association in which tables with single small entries are up-weighted but those with skewed margins are down-weighted according to the relative entropy 
among the tables of the same odds-ratio. 
\\ \\ \\ \\
{\bf Keywords:} two by two probability tables, measures of association, entropy 

\newpage
\section{Introduction}
2x2 tables of binary markers with random margins are intriguing in several respects: 
First, there is a confusing plethora of measures of association in 2x2 tables with random margins that are used in statistical practice. Their relative merit is unclear. 
Some of them were developed for 2x2 tables with fixed margins and then extended to the case considered here. Measures typically agree in the ordering by strength of 
association on 2x2 tables that have diagonal symmetry and in case of independence. But they markedly differ in asymmetric tables or in tables which are "far from independence". 
We develop a unified framework to analyse, visualise and compare measures of association in 2x2 probability tables especially with respect to their dependence on the
margins.

Second, 2x2 tables "far from independence" may approximate logical forms like logical  equivalence (one diagonal is zero) or implication (one entry zero). The task of 
selecting particularly interesting and informative tables among a large number of tables  is often encountered in the analysis of data consisting of high 
dimensional binary patterns (e.g. linkage disequilibrium of SNPs, patterns of aberration at various DNA loci, patterns of protein expression etc.). We suggest a principled 
approach for picking tables which approximate logical relations. This approach relies on an entropy-based weighting of tables and aims to improve existing measures often used in Genetical Statistics. 

Defining and justifying measures and estimating them from empirical data are radically different tasks. We have investigated methods of estimating measures of association in a separate paper 
\citep{SchHas}. Here we deal exclusively with abstract 2x2 probability models and their mathematical structure.  

\section{Mathematical structure of 2x2 probability models}

2x2 tables of binary markers with random margins can be considered as tetranomial distributions with a symmetry structure. Symmetry of 2x2 tables can be described by the dihedral 
group $D_{4}$ generated by the transposition of the binary markers (matrix transposition) and transposition of their values (transposition of columns or rows).

We consider the manifold ${\mathbb T}$ of all non-degenerate tetranominal probability models which we write in two by two lay-out: ${\mathbb T}$ consists of all two by 
two matrices $t$ with entries $p_{ij}\in {\mathbb R}$, ($i,j\in \{0,1\}$) subject to the constraints
$p_{ij}>0$, $\sum_{i,j}p_{ij}=1$. The $p_{ij}$ denote the probabilities of the corresponding combination of the states of two binary markers $i$ and $j$. In the following, 
we abbreviate $\sum_{i=0}^1\sum_{j=0}^1=\sum_{i,j}$, $p_{i.}=p_{i0}+p_{i1}$ and
$p_{.j}=p_{0j}+p_{1j}$. The margins $p_{i.}$ and $p_{.j}$ give the marginal distributions of the marker $i$ and $j$ respectively. 

In ${\mathbb T}$ we have several relevant submanifolds. There is a marked point $m_{0}$, namely the midpoint  $\left({1/4\atop 1/4}{1/4\atop 1/4}\right)$. There is the 1-dimensional submanifold 
${\mathbb D\mathbb S}$ of all tables with diagonal symmetry of the form $\left({a\atop b}{b\atop a}\right)$. And there is the 2-dimensional submanifold ${\mathbb I \mathbb N \mathbb D}$ of 
independent tables with $p_{ij}=p_{i.}\cdot p_{.j}\quad\forall i,j$.

By $\overline{{\mathbb T}}$ we denote the closure of ${\mathbb T}$. The border $\partial\overline{{\mathbb T}}= \overline{{\mathbb T}} - \mathbb T $ consists of tables with at least one zero: 
four two dimensional sides $\left\{{p_{ij}=0}\right\}$ for any $i,j$, six one dimensional edges of vanishing rows $\left\{p_{.j}=0\right\}$, vanishing columns $\left\{p_{i.}=0\right\}$ and two 
vanishing diagonals $\left\{p_{00}=p_{11}=0\right\}$, $\left\{p_{01}=p_{10}=0\right\}$ as well as four triple zero vertices $\left\{p_{ij}=1\right\}$. 

Manipulating the margins defines an additional structure on ${\mathbb T}$. We can multiply rows or columns with positive numbers and renormalise: Formally, consider the group 
$G=(\mathbb{R^+}\times\mathbb{R^+},\cdot)$ with component-wise multiplication. \\
For every $(\mu,\nu) \in \mathbb{R^+}\times\mathbb{R^+}$ we define a map: $g(\mu,\nu):\mathbb{T}\longrightarrow \mathbb{T}$
\begin{eqnarray} 
t=\left({p_{00}\atop p_{10}}{p_{01}\atop p_{11}}\right) \longmapsto g(\mu,\nu)(t)&=&\frac{1}{\mu\nu p_{00}+\mu p_{01}+\nu p_{10}+p_{11}}\left(\begin{array}{ll} \mu\nu p_{00} & \mu p_{01} \\
\nu p_{10} &  p_{11} \end{array}\right) 
\end{eqnarray}
Since $g(\mu,\nu) \circ g(\mu^\prime,\nu^\prime) = g(\mu\cdot\mu^\prime,\nu\cdot\nu^\prime)$ and $g(1,1)={}$Id$_\mathbb{T}$ this defines a G-group action on $\mathbb{T}$. 

Lying in the same group orbit defines an equivalence relation on $\mathbb{T}$: We say  two elements $t_1,t_2 \in \mathbb{T}$ are equivalent $t_1 \sim t_2$ if and only if there are 
$(\mu,\nu)\in\mathbb{R^+}\times\mathbb{R^+}$ with $g(\mu,\nu)(t_1)=t_2$. G-Orbits are diffeomorph to $\mathbb{R^+}\times\mathbb{R^+}$. 

A real function $\eta: \mathbb{T} \rightarrow \mathbb{R}$ is $G$-invariant if $\eta(t)=\eta(g(\mu,\nu)(t))$ for all $(\mu,\nu) \in \mathbb{R^+}\times\mathbb{R^+}$.
\\ \\ \noindent
{\bf Proposition 1 (odds-ratio):}\\
{\it 
a) The odds-ratio $\lambda:\mathbb{T} \rightarrow \mathbb{R}; \ t=\left({p_{00}\atop p_{10}}{p_{01}\atop p_{11}}\right) \mapsto \lambda(t)=\frac{p_{00}p_{11}}{p_{01}p_{10}}$ is G-invariant.\\ 
b) The odds-ratio classifies the G-orbits. Let $\tilde{\mathbb{T}}$ be the quotient space of $\mathbb{T}$ by the equivalence relation induced by G. $\lambda$ induces a bijective map  
$\tilde{\lambda}: \tilde{\mathbb{T}}\stackrel{~}{\rightarrow}\mathbb{R^+}$.\\
c) The inverse mapping $\tilde{\lambda}^{-1}: \mathbb{R^+} \rightarrow \tilde{\mathbb{T}}$ can be described by 
$l  \longmapsto \left[\left(\begin{array}{ll} \frac{\sqrt{l}}{2\cdot(1+\sqrt{l})} & \frac{1}{2\cdot(1+\sqrt{l})} \\
\frac{1}{2\cdot(1+\sqrt{l})} &  \frac{\sqrt{l}}{2\cdot(1+\sqrt{\l})} \end{array}\right)\right]$\\
d) Every G-invariant function $\eta: \mathbb{T} \rightarrow \mathbb{R}$ can be written as a function of $\lambda$, namely $\eta=(\tilde{\eta}\circ\tilde{\lambda}^{-1})\circ\lambda$. }
\\ \\ \noindent
{\it Proof:} a) is easily verified. b) Every equivalence class [t] in $\tilde{\mathbb{T}}$ has a representant with margins $\frac{1}{2}$, namely 
$[g(\sqrt{\frac{p_{10}p_{11}}{p_{00}p_{01}}},\sqrt{\frac{p_{01}p_{11}}{p_{00}p_{10}}})(t)]$ which has the form given in c). d) is trivial. \hfill $\Box$
\\ \\
We next define new coordinates on ${\mathbb T}$ to make use of this insight. 
\\\\
{\bf Proposition 2 (Margin transformation coordinates on ${\mathbb T}$ highlighting the G-action and its invariant):} {\it The map $\Theta: \mathbb T\rightarrow \mathbb{R}^{3}$
\begin{eqnarray} 
t=\left({p_{00}\atop p_{10}}{p_{01}\atop p_{11}}\right) \longmapsto \Theta(t)=(\ln\sqrt{\frac{p_{00}p_{11}}{p_{01}p_{10}}},\ln\sqrt{\frac{p_{00}p_{01}}{p_{10}p_{11}}},\ln\sqrt{\frac{p_{00}p_{10}}{p_{01}p_{11}}})
\end{eqnarray} 
is a diffeomorphism. \\
The inverse $\Psi=\Theta^{-1}: \mathbb{R}^{3} \rightarrow \mathbb T$ is given by}
\begin{eqnarray*}
\Psi(x,y,z) &=&  g\left(e^y,e^z\right)) \left ( \left(\begin{array}{ll} \frac{e^x}{2\left(1+e^x\right)} & \frac{1}{2\left(1+e^x\right)} \\
\frac{1}{2\left(1+e^x\right)} &  \frac{e^x}{2\left(1+e^x\right)} \end{array}\right)\right) \\
&=&\frac{1}{e^{x+y+z} + e^x + e^y + e^z} 
\left(\begin{array}{ll} e^{x+y+z} & e^y \\ e^z & e^x \end{array}\right) 
\end{eqnarray*}
\\ \noindent
In these new coordinates, $x$ corresponds to the log odds-ratio, while $y$ and $z$ determine the G-transformation that maps the table to diagonal symmetry.
In addition, the midpoint $m_0$ corresponds to the origin $(0,0,0)$. G-orbits (odds-ratio = constant) correspond to planes 
${\left\{a\right\}}\times\mathbb{R}^{2}$. In particular, the submanifold of independent tables ${\mathbb I \mathbb N \mathbb D}$ maps to 
${\left\{0\right\}}\times\mathbb{R}^{2}$. The tables with diagonal symmetry ${\mathbb D\mathbb S}$ form the line  
$\mathbb{R}\times{\left\{0\right\}}\times{\left\{0\right\}}$. Transposing rows and columns of a table is equivalent to transformations $y \to -y$
and $z \to -z$, while matrix transposition is equivalent to the transformation $y \leftrightarrow z$. 

Let $\bar{\mathbb{R}}:= \mathbb{R} \cup \{-\infty,+\infty\} $ be the two point compactification of $\mathbb{R}$. $\bar{\mathbb{R}}^3$ is a compactification of $\mathbb{R}^3$ as a cube. We use a short hand notation to describe the boundaries abbreviating $+\infty$ as ''+'', $-\infty$ as ''-'' and any finite real number as ''*''. The eight vertices  $V=\left\{(\pm\pm\pm)\right\}$ split into two sets of four: $V_g=\{(+++),(+--),(-+-),(--+)\}$ and $V_b=\{(---),(-++),(+-+),(++-)\}$. 
\\\\
{\bf Proposition 3 (Extension to the borders):}
{\it $\Psi$ and $\Theta$ considered as set valued functions can be extended to $\bar{\mathbb{R}}^3$ respectively $ \overline{{\mathbb T}}$.  They remain inverse to each other. The mappings of the borders can be characterized as follows:
\begin{itemize}
\item The vertices $V_g$ together with their respective adjacent edges map to the vertices in $\overline{{\mathbb T}}$.
\item The faces of $\overline{{\mathbb T}}$ correspond to the vertices $V_b$.
\item The faces $\left(\pm**\right)$  of the cube map to the diagonal edges $p_{00}=p_{11}=0$ and $p_{01}=p_{10}=0$ in $\overline{{\mathbb T}}$.
\item The faces $\left(*\pm*\right)$, $\left(**\pm\right)$  correspond to tables with vanishing rows $\left\{p_{.j}=0\right\}$ or vanishing columns $\left\{p_{i.}=0\right\}$ in $\overline{{\mathbb T}}$
respectively. 
\end{itemize}}
This behaviour is illustrated in figure 1.
These different compactifications will later be used to characterise the limit behaviour of association measures. It will turn out that the limit behaviour can be easier described using the margin transformation coordinates.  

\section{Measures of association}
We will now investigate various measures of associations between two binary markers. First we define the objects of interest.\\ \\ \noindent
{\bf Definition (Measures of association):}
{\it A measure of association between binary markers is a continuous function $\eta: \mathbb{T} \rightarrow \mathbb{R}$ with the following properties: \\
a) $\eta$ is zero on independent tables.\\
b) $\eta$ is a strictly increasing function of the odds-ratio when restricted to tables with fixed margins.\\
c) $\eta$ respects the symmetry group $D_{4}$, namely:\\
c1) $\eta$ is symmetric in the markers, i.e. invariant to matrix transposition.\\
c2) $\eta$ changes sign when states of a marker are transposed (row or column transposition). }\\
\\ \noindent A measure of association is \textit{standardised} if its range is restricted to $(-1,1)$.

\subsection{Measures based on the odds-ratio}
The odds-ratio {\it Odds-ratio $\lambda$} \citep{Ed}:
\begin{eqnarray*}
\lambda&=&\frac{p_{00}p_{11}}{p_{01}p_{10}}
\end{eqnarray*}
can be used to define measures of association. As $\lambda$ is G-invariant, monotone transformations automatically fulfill condition b) of the definition.
\\ \noindent
Standardised measures of association derived from the odds-ratio include 
{\it Yule's $Q$} \citep{Yule}:
\begin{eqnarray*}
Q&=&\frac{\lambda-1}{\lambda+1}
\end{eqnarray*}
and 
{\it Yule's $Y$} \citep{Yule}:
\begin{eqnarray*}
Y&=&\frac{\sqrt{\lambda}-1}{\sqrt{\lambda}+1} 
\end{eqnarray*}
Obviously, both $Q$ and $Y$ are measures of association. 
Similar to the odds-ratio, both are extremal if one of the $p_{ij}$ tends to zero.

\subsection{Measures based on additive deviations from independence given the margins}
Fixing margins results is a one dimensional submanifold of tables that can be additively parametrised by a parameter D. \\
All such tables have the form: 
\begin{eqnarray*}
\left(\begin{array}{ll} {p_{0.}\cdot p_{.0}+D} & {p_{0.}\cdot p_{.1}-D} \\
{p_{1.}\cdot p_{.0}-D} & {p_{1.}\cdot p_{.1}+D}\end{array}\right)
\end{eqnarray*}
$D=p_{00}p_{11}-p_{01}p_{10}\gl p_{00}-p_{0.}p_{.0}$ describes the additive deviation from the independent table with the given margins. 
This measure is zero in case of independence of the markers but extremal values depend on the margins.\\
\\
{\it Lewontin's $D^\prime$} \citep{Lew}:
The measure $D^\prime$ is a standardisation of the original measure $D$:
\begin{eqnarray*}
D^\prime&=&\frac{D}{D_{max}} \qquad\mbox{where}\quad D_{max}=\left\{\begin{array}{l@{\quad\mbox{if}\quad}l}
\min\left\{p_{0.}p_{.1},p_{.0}p_{1.}\right\} & D\ge 0 \\ 
\min\left\{p_{0.}p_{.0},p_{1.}p_{.1}\right\} & D<0 \end{array} \right.
\end{eqnarray*}
Lewontin's $D^\prime$ ranges from $-1$ to $1$ and tends to these values if at least one of the $p_{ij}$ tends to zero. \\
$D^\prime$ is widely used in genetics to measure linkage disequilibrium. When a new SNP emerges in a population by a single mutation event, the new allele is exclusively found in 
conjunction with only one of the two alleles of already existing SNPs. As long as no recombination events occurs, the new SNP remains in complete linkage disequilibrium with the 
other SNPs. The corresponding 2x2 tables feature a single zero cell. Thus in this context a measure is needed that is extremal whenever a single entry tends to zero.

Since $D_{max}$ is constant for tables with fixed margins and $D$ increases with increasing odds-ratio, $D^\prime$ is a monotone function of the odds-ratio for constant margins. Symmetry 
is obvious. Hence, $D^\prime$ is a standardised measure of association.
\\
\\
{\it Correlation coefficient $r$} \citep{HR}: The correlation 
coefficient applied to binary data has similar popularity in genetics as $D^\prime$. It ranges also from $-1$ to $1$, but, in contrast to $D^\prime$, the absolute value $1$ is obtained when a diagonal of $t$ tends to zero:
\begin{eqnarray*}
r&=&\frac{D}{\sqrt{p_{0.}p_{.0}p_{1.}p_{.1}}} \gl \frac{p_{00}p_{11}-p_{01}p_{10}}{\sqrt{p_{0.}p_{.0}p_{1.}p_{.1}}}
\end{eqnarray*}
With reasoning similar as for $D^\prime$, $r$ is a standardised measure of association.
\\
\\
\noindent
{\bf Proposition 4 (Equality of $r$, $D^\prime$ and $Y$ on diagonal tables): } {\it The measures $r$, $D^\prime$ and $Y$ coincide on the set of diagonal tables, i.e. tables with pair-wise 
equal diagonal elements.} \\ \\ \noindent
{\it Proof:} This follows directly after calculating these measures for the tables $t=\frac{1}{2a+2b}\left( {a\atop b}{b\atop a} \right)$, $a,b>0$. \hfill $\Box$

\subsection{Measures based on information theory}
The {\it mutual information} \citep{Shan} is defined as the difference between the information of the given table and the independent table with the same margins.
\begin{eqnarray*}
\MutInf&=&\sum_{i,j}p_{ij} \cdot \log_2(p_{ij}) - \sum_{i}p_{i.} \cdot \log_2(p_{i.})-\sum_{j}p_{.j} \cdot \log_2(p_{.j})
\end{eqnarray*}
$\MutInf$ takes values only in $[0,1]$. In order to make it a measure of association according to our definition, we define a signed version:  
\begin{eqnarray*}
\sMutInf&=&\sign(D) \cdot \MutInf
\end{eqnarray*}
\\ \noindent
{\bf Proposition 5:} {\it $\sMutInf$ is a standardised measure of association. }
\\ \\
\noindent
{\it Proof:} The symmetry of this measure is clear. To show that $\sMutInf$ is a monotone function of the odds-ratio, we consider the tables 
$t_\varepsilon=\left( {p_{00}+\varepsilon \atop p_{10}-\varepsilon} {p_{01}-\varepsilon \atop p_{11}+\varepsilon} \right)$ for a sufficiently small $\varepsilon>0$. 
These tables have the 
same margins as the table $t=\left( {p_{00} \atop p_{10}} {p_{01} \atop p_{11}} \right)$ but higher odds-ratios. Assume that $\lambda>1$, we see that
$\left.\frac{d}{d\varepsilon}\right|_{\varepsilon=0}\sMutInf(t_\varepsilon)=\log_2 \lambda >0$. Hence $\sMutInf$ is monotone, and thus, a measure of association. \hfill $\Box$
\\ 
\\
\noindent
$\MutInf$ approaches 1 only if approaching $\left({1/2\atop 0}{0\atop 1/2}\right)$ while $r$ approaches 1 by approaching tables of the 
form $\left({a\atop 0}{0\atop b}\right)$, $a,b>0$. 

\subsection{Counter example}
{\it Kappa coefficient} \citep{Coh}: The Kappa coefficient which is useful in quantifying the agreement between two raters is defined as: 
\begin{eqnarray*}
\kappa&=& \frac{p_{00}+p_{11}-p_{0.}p_{.0}-p_{1.}p_{.1}}{1-p_{0.}p_{.0}-p_{1.}p_{.1}}
\end{eqnarray*}
Kappa is {\bf not} a measure of association. Although it fulfils the condition of monotonicity, it is not symmetric. 

\section{Comparing measures of association}

We use the coordinates introduced in Proposition 2 in order to describe and visualise how measures of association depend on the margins. In particular we study measures of association $\eta$ restricted to $x$=const i.e. for fixed odds-ratios. The restricted functions will be denoted $\eta_x$ and called margin weighting functions. We characterise the shape of the margin weighting functions and study their limiting behaviours and extensibility to the compactification $\bar{\mathbb{R}}^3$ in comparison to $ \overline{{\mathbb T}}$.   
\\
The association measure $r$ expressed in margin transformation coordinates reads:
\begin{eqnarray}
r\left(x,y,z\right)&=&\frac{\left(e^{2x}-1\right)e^{y+z}}
{\sqrt{\left(e^{x+y+z}+e^y\right)\left(e^{x+y+z}+e^z\right)
\left(e^x+e^y\right)\left(e^x+e^z\right)}}
\end{eqnarray}
\noindent
The margin weighting function of $r$ for odds-ratio $\lambda=40$ is shown in figure 2. 
\\
\\
\noindent
{\bf Proposition 6 (Margin weighting function for $r$):} {\it
For all $x \in \mathbb{R}\setminus \left\{0 \right\}$:\\
a) ${r}_{x}$ has exactly one extremum at the origin $\left(y,z\right)=\left(0,0\right)$, corresponding to the diagonal symmetric table with the fixed odds-ratio.\\
b) $\lim_{\left\|(y,z)\right\| \rightarrow \infty }r_x=0$. \\
c) $\lim_{x\to\pm\infty} r_x=\pm 1$ \\
d) $r$ can be extended to $\bar{\mathbb{R}}^3$ except for the lines $(\pm,\pm,*)$ and $(\pm,*,\pm)$ and the vertices $V$. \\
e) $r$ can be extended to $ \overline{{\mathbb T}}$ except for the vertices. }
\\ \\ \noindent
{\it Proof:} see appendix.
\\ \\ \noindent
The measure $r$ down-weights tables with skewed margins.
\\ \\
The association measure $D^\prime$ expressed in margin transformation coordinates reads:
\begin{eqnarray}
D^\prime\left(x,y,z\right)&=&\frac{\left(e^{2x}-1\right)e^{y+z}}{D_{max}}\qquad\mbox{where} \\
D_{max}&=&\left\{\begin{array}{l@{\quad:\quad}l}
\left(e^{x+y+z}+e^y\right)\left(e^x+e^y\right) & x>0,\quad y<z \\
\left(e^{x+y+z}+e^z\right)\left(e^x+e^z\right) & x>0,\quad y\ge z \\
\left(e^{x+y+z}+e^y\right)\left(e^{x+y+z}+e^z\right) & x<0,\quad y<-z \\
\left(e^x+e^y\right)\left(e^x+e^z\right) & x<0,\quad y\ge -z \\
\end{array} \right. \nonumber
\end{eqnarray}
The margin weighting function of $D^\prime$ for odds-ratio $\lambda=40$ is shown in figure 3. 
\\
\\
\noindent
{\bf Proposition 7 (Margin weighting function for $D^\prime$):}  {\it 
For all $x \in \mathbb{R}\setminus \left\{0 \right\}$:\\
a) ${D^\prime}_{x}$ has a non-differentiable edge along the diagonal $y=z$ for $D^\prime>0$ and along the diagonal $y=-z$ for $D^\prime<0$. There is a non-smooth saddle point in the origin.\\
b) 
\begin{eqnarray*}
\lim_{y\to\pm\infty}D^\prime_x&=&\left(e^{2x}-1\right)\cdot\left\{\begin{array}{l@{\quad:\quad}l}
\left(e^{2x}+e^{x\pm z}\right)^{-1} & x>0 \\
\left(e^{x\mp z}+1\right)^{-1} & x<0 
\end{array}\right. \\ 
\lim_{z\to\pm\infty}D^\prime_x&=&\left(e^{2x}-1\right)\cdot\left\{\begin{array}{l@{\quad:\quad}l}
\left(e^{2x}+e^{x\pm y}\right)^{-1} & x>0 \\
\left(e^{x\mp y}+1\right)^{-1} & x<0 
\end{array}\right. 
\end{eqnarray*}
Thus, limit functions have a range of $\left(0,1-e^{-2x}\right)$ for $x>0$ and 
$\left(e^{2x}-1,0\right)$ for $x<0$, where 0 is obtained for 
$y\to\pm\infty$, $z\to\pm\infty$, $x>0$ and $y\to\mp\infty$, $z\to\pm\infty$, $x<0$.\\
c) $\lim_{x\to\pm\infty} D^\prime_x=\pm 1$ \\
d) $D^\prime$ can be extended to $\bar{\mathbb{R}}^3$ except for the vertices $V_g$. \\
e) $D^\prime$ can be extended to $ \overline{{\mathbb T}}$ except for the edges and vertices.}
\\ \\ \noindent
{\it Proof:} see appendix.
\\ \\ \noindent
$D^\prime$ gives higher weights to certain tables without diagonal symmetry. The measure 
up-weights or down-weights tables with skewed margins depending on the position of zeros which occur in the limiting tables (see figure 3). Comparing d) and e) one recognizes that the introduction
of the odds-ratio as coordinate allows extending $D^\prime$ to limit tables with vanishing colums or rows.
\\ \\
The association measure $\sMutInf$ can also be written in margin transformation coordinates but this is skipped due to 
the lengthy formula.
The margin weighting function of $\sMutInf$ for odds-ratio $\lambda=40$ is shown in figure 4. \\
\\
\noindent
{\bf Proposition 8 (Margin weighting function for $\sMutInf$):}   {\it 
For all $x \in \mathbb{R}\setminus \left\{0 \right\}$:\\
a) $\sMutInf_{x}$ has exactly one maximum at the origin $\left(y,z\right)=\left(0,0\right)$.\\
b) $\lim_{\left\|(y,z)\right\| \rightarrow \infty }\sMutInf_x=0$. \\
c) $\lim_{x\to\pm\infty} \sMutInf_x=\pm\left(\log_2\left( e^{y\pm z}+1\right)-\frac{e^{y\pm z}}{e^{y\pm z}+1}\log_2 e^{y \pm z}  \right)$\\ 
Thus $\sMutInf_{x}\to \pm 1$ for $y=\mp z$ and $x \to\pm\infty$ respectively. \\
d) $\sMutInf$ can be extended to $\bar{\mathbb{R}}^3$ except for the vertices $V_b$. \\
e) $\sMutInf$ can be extended completely to $ \overline{{\mathbb T}}$. }
\\ \\ \noindent
{\it Proof:} see appendix.
\\ \\ \noindent
Thus, similarly to $r$, $\sMutInf$ down-weights tables with skewed margins (see figure 4). 
\\ \\ \noindent
The association measure $Y$ in margin transformation coordinates can be simply written as:
\begin{eqnarray}
Y\left(x,y,z\right) &=& \tanh \frac{x}{2}
\end{eqnarray}
\noindent
{\bf Proposition 9 (Margin weighting function for $Y$):} {\it 
For all $x \in \mathbb{R}$:\\
a) $Y_x$ is constant.\\
b) $\lim_{\left\|(y,z)\right\| \rightarrow \infty }Y_x=\tanh \frac{x}{2}$. \\
c) $\lim_{x\to\pm\infty} Y_x=\pm 1$\\ 
d) $Y$ can be extended completely to $\bar{\mathbb{R}}^3$. \\
e) $Y$ can be extended to $\overline{{\mathbb T}}$ except for edges and vertices corresponding to vanishing rows or columns. }
\\ \\ \noindent
{\it Proof:} is trivial. \hfill$\Box$

\section{Entropy}

Among tables of a fixed odds-ratio, we look for a principled approach to prefer interesting tables and down-weight obscure "junk" tables. As a candidate we study the table entropy on $\mathbb{T}$. 
The entropy function $H: \mathbb{T} \rightarrow \mathbb{R}$ is defined as the negative expectation of the loglikelihood of the tables:
\begin{eqnarray*}
H\left(\left({p_{00}\atop p_{10}}{p_{01}\atop p_{11}}\right)\right)&:=& - \left(p_{00} \cdot \log_2(p_{00}) +p_{01} \cdot \log_2(p_{01}) +p_{10} \cdot \log_2(p_{10}) +p_{11} \cdot \log_2(p_{11})\right)
\end{eqnarray*}
Why is entropy a candidate to select among tables? It can be charcterised in multiple ways:
For general finite discrete distributions the entropy was introduced by Shannon (1948) \citep{Sh}. Shannon characterised $H$ by a set of postulates to measure the uncertainty in a discrete distribution: 
\\
\\ 
\noindent
{\bf Shannon's characterisation of Entropy:} {\it 
If functions $H_n(p_1,...,p_n)$ with $p_i \geq 0, \sum p_i =1, n\geq 2$ satisfy the conditions\\
a) $H_2(p,1-p)$ is a continuous positive function of $p$.\\
b) $H_n(p_1,...,p_n)$ is symmetric, i.e. invariant under permutations of the $p_1,...,p_n$ for all n.\\
c) $H_n(p_1,...,p_n)=H_{n-1}(p_1+p_2,p_3,...,p_n)+(p_1+p_2)\cdot H_2(\frac{p_1}{(p_1+p_2)},\frac{p_2}{(p_1+p_2)}  )$\\
then $H_n(p_1,...,p_n)=-K\cdot \sum p_i \log_2(p_i)$ for some $K>0$. }
\\
Tables with high entropy are interesting as they have high uncertainty and "surprise value".
\\ \\ \noindent
Jaynes \citep{Jaynes} gives an independent combinatorial characterisation: When we sample sequentially from a table $t \in \mathbb{T}$ we obtain a vector of observations of length $N$, 
which we summarise as a frequency table $\hat{t}_N=1/N\cdot\left({n_{00}\atop n_{10}}{n_{01}\atop n_{11}}\right)$ . Each frequency table $\hat{t}_N$ is characterised by the number 
$W(\hat{t}_N)= \frac{N!}{n_{00}!n_{01}!n_{10}!n_{11}!}$ of sequences which realise $\hat{t}_N$.  Intuitively, tables that can be realised in multiple ways are more plausible than those that can 
be realised only by few sequences. 
We can use Stirlings formula for $n!$ to approximate $W(\hat{t}_N)$. In the limit $N \rightarrow \infty$,   $\hat{t}_N \rightarrow t$ in probability and $1/N\cdot\log(W(t^N)) \to H(t)$. Thus the 
entropy describes the combinatorial plausibility of a table.

Given a set of distributions fulfilling certain constraints, Jaynes \citep{Jaynes} proposes to pick the corresponding maximum entropy distribution as the most uncommitted and prototypical 
distribution. 
Looking at the margin weighting function of the entropy leads to a surprise:

Recall that Lambert's W-function is defined as the inverse function to $x \exp x$. $W$ is a multi-branch function since $y=x \exp(x)$ has two solutions for $y \in (-1/e,0)$. We can prove the following:
\\ \\
\noindent 
 {\bf Theorem 1 (magic odds-ratio):} {\it 
Define the "magic odds-ratio" by $L_{magic}=W(1/e)^{-2} \approx 12.89 $. 
Let $L > 1$. The entropy $H$ restricted to the submanifold of constant odds-ratio $L$ in $\mathbb{T}$
\begin{itemize}
	\item has a single maximum at the diagonal table of odds-ratio $L$ if $1 < L \le L_{magic}$.

	\item has a saddle point at the diagonal table of odds-ratio $L$ and two "L-shaped" tables as maxima which transpose with matrix transposition if $L_{magic} < L$.
\end{itemize}
\noindent "L-shaped" means that for $L \rightarrow \infty$ one of the maxima approaches the table $\left( {1/3\atop 0}{1/3\atop 1/3}\right)$.  
For the case $L<1$ a similar result can be derived by transposing principal and secondary diagonals. }
\\ \\ \noindent
{\it Proof:}  There are two constraints to be considered, one of them not linear in $p_{ij}$:
\begin{eqnarray}
 \ln(p_{00}) - \ln(p_{01}) - \ln(p_{10}) +\ln(p_{11}) &=& \ln(L) \\
 p_{00} + p_{01} + p_{10} + p_{11} &=& 1 
\end{eqnarray}
Using Langrange multipliers, critical tables of H restricted to odds-ratio equals $L$ can be expressed in terms of Lambert's W function. The bifurcation occurs for $L_{magic} < L$ because Lambert's W is multi\-branched. See appendix for details.
\\ \\
\noindent 
This theorem suggests that the "magic odds-ratio" is a natural cutpoint between weak and strong association. For weak association $L < L_{magic}$, interesting tables are those near ${\mathbb D\mathbb S}$. For strong association $L_{magic} < L$, particularly interesting tables are those that approach "L-shape", i.e. those in which one cell differs in magnitude from the three others. 

\section{An entropy-based measure of association}

Using these insights on the entropy of a table, in this section we aim to define a measure of association with similar properties to $D^\prime$, $Y$ but better limit behaviour, i.e. 
the measure should down-weight tables with almost vanishing rows or columns or single entries. These tables are denoted as {\it junk tables} in the following. We have seen in the 
last sections that $D^\prime$ and $Y$ could be large for these tables. 

We also like to recall that both, $D^\prime$ and $Y$ become extremal if the table features a single entry equals zero while $r$, $\sMutInf$ require a vanishing diagonal. 
We like to retain this property for a new measure to be defined. Another feature to be retained is the agreement of measures for diagonal tables which holds for $Y$, $D^\prime$ 
and $r$. 

According to our definition, an important property of a measure of association is that it is a monotone function of the odds-ratio when the margins are kept fixed. 
For the entropy, one can prove the following lemma:
\\ \\ \noindent
{\bf Lemma 1 (Monotony of the entropy difference):} {\it Let $H$ be the entropy of $t$ and $H_{diag}$ be the entropy of the corresponding diagonal table of the same odds-ratio $\lambda$. 
Then, $H_{diag}-H$ is monotonically decreasing for increasing $\lambda>1$ and constant margins. }
\\ \\ \noindent
{\it Proof:} see appendix. 
\\ \\ 
\noindent
As a direct consequence of this lemma, it is easy to see that:
\\ \\
\noindent
{\bf Corollary:} 
\begin{eqnarray}
\HS_n &:=& \sign Y \left|Y\right|^{\exp n\left(H_{diag}-H\right)}
\end{eqnarray}
{\it is a measure of association for arbitrary $n\ge0$.}
\\ \\ \noindent
This newly defined measure fulfils all above mentioned properties: It coincides with $Y$, $D^\prime$, $r$ at diagonal tables, is extremal for tables with a single zero, up-weights 
L-shaped tables for large odds-ratios in the sense that $\HS_n>Y$ and down-weights junk-tables in the sense that $\HS_n<Y$ at the margins (proof see below). However, the down-weighting is 
imperfect as $\HS_n>0$ for junk-tables.

The parameter $n$ can be chosen in order to define the degree of up- and down-weighting. According to our observations, $n=4$ is a reasonable choice resulting in 
a satisfactory down-weighting of junk tables (see later). 

The measure $\HS_n$ can be written in margin transformation coordinates using 
\begin{eqnarray*}
H_{diag}\left(x,y,z\right)&=&1+\log_2\left(1+e^x\right)-\frac{x}{\ln 2\left(e^{-x}+1\right)}\\
H\left(x,y,z\right)&=&\log_2\left(e^{x+y+z}+e^x+e^y+e^z\right)-\frac{\left(x+y+z\right) e^{x+y+z} + xe^x + ye^y+ ze^z}{\ln 2\left(e^{x+y+z}+e^x+e^y+e^z\right)}
\end{eqnarray*}
At figure 5 we present the margin weighting functions of $\HS_n$ for $\lambda=5$ and $\lambda=40$. 
These functions can be easily characterised using the results of the previous section:
\\ \\ \noindent
{\bf Proposition 10 (Margin weighting function for $\HS_n$):} {\it
For all $x \in \mathbb{R}\setminus \left\{0 \right\}$:\\
a) For $x\in\left(-1-W\left(1/e\right),1+W\left(1/e\right)\right)$, 
${\HS_n}_{x}$ has exactly one maximum at the origin $\left(y,z\right)=\left(0,0\right)$. If $x<-1-W\left(1/e\right)$ or 
$x>1+W\left(1/e\right)$, ${\HS_n}_{x}$ has a saddle-point at the origin and two extrema elsewhere. At these extrema,
the elements of one diagonal are equal while at the other diagonal there is one (small) element. \\ \noindent
b) ${\HS_n}_{x}$ has the following limit functions
\begin{eqnarray*}
\lim_{\left\|(y,z)\right\| \rightarrow \infty}  {\HS_n}_x &=& \sign\left( \tanh \frac{x}{2}\right) \left|\tanh 
\frac{x}{2}\right|^{\exp n 
\left\{1+\log_2\left(1+e^x\right)-\frac{x}{\ln 2\left(e^{-x}+1\right)}+p\log_2 p + \left(1-p\right)\log_2 \left(1-p\right)\right\}}
\end{eqnarray*}
where $p=\left(1+e^{x\pm z}\right))^{-1}$ for $y\to\pm\infty$ or $p=\left(1+e^{x\pm y})\right)^{-1}$ for 
$z\to\pm\infty$ respectively. Thus, the limit functions have an extremum at $p=0.5$ that is $z=\mp x$ for $y\to\pm\infty$ and $y=\mp x$ for $z\to\pm\infty$ respectively. \\
c) $\lim_{x\to\pm\infty}{\HS_n}_x =\pm 1$ \\
d) ${\HS_n}_{x}<Y_x$ at the margins, i.e. $\HS_n$ down-weights junk-tables. \\
e) $\HS_n$ can be extended completely to $\bar{\mathbb{R}}^3$. \\
f) $\HS_n$ can be extended to $\overline{{\mathbb T}}$ except for edges and vertices corresponding to vanishing rows or columns. \\
g) For all $x \in \mathbb{R}$, $\HS_n$ coincides with $Y$, $D^\prime$, $r$ at diagonal tables. }
\\ \\ \noindent
{\it Proof:} a) follows from the Theorem 1. b) is easy to see taking the limit of the tables first. c) is clear 
since $\lim_{x\to\pm\infty} \tanh\frac{x}{2}=\pm 1$ and the exponent is finite.
d) holds since $H_{diag}>1$ and $H\le 1$ at the margins of finite $x$. 
e) and f) are consequences of b) and c). g) is obvious. \hfill $\Box$

\section{Examples of tables and corresponding association measures}
We now study the behaviour of the measures $Y$, $r$, $D^\prime$ and the newly proposed measure $\HS_4$ for a variety of selected tables (see table 1). For this purpose, we study the 
odds-ratios $\lambda\in\left\{1,2,5,10,20,50,100\right\}$
and consider the following tables for $x=\ln\sqrt{\lambda}$: 
\begin{itemize}
\item The diagonal table ($y=z=0$).
\item An L-shaped table, characterized by $y=x$, $z=-x$.
\item A junk table with $y=10$, $z=-y$ corresponding to $p_{01}\approx 1$.
\item A junk table with $y=10$, $z=-x$ corresponding to $p_{00}\approx p_{01}\approx 0.5$.  
\item A junk table with $y=10$, $z=y$ corresponding to $p_{00}\approx 1$. 
\end{itemize}
We also like to remark that the table with three equal entries has maximum entropy if $\lambda\to\infty$. 

Per definition of a measure, for $\lambda=1$ all measures equals zero independent of the concrete realization of the table. Since $Y$ is 
based on the odds-ratio. $Y$ is constant for all tables of the same odds-ratio. $Y$, $r$, $D^\prime$ and $\HS_4$ always coincide at diagonal tables.
$r$ is maximal at diagonal tables and becomes small for all kinds of junk tables. $D^\prime$ is always greater for L-shaped tables  
than for diagonal tables. $D^\prime$ is close to zero in case of $p_{00}\approx 1$ but could become large for $p_{01}\approx 1$ which is highly counter-intuitive. 
$\HS_4$ also becomes larger for L-shaped tables compared to diagonal tables if $\lambda$ is large. In contrast to $D^\prime$, $\HS_4$ is close to zero 
for both junk configurations $p_{00}\approx 1$ and $p_{01}\approx 1$ respectively. The limit tables have a maximum of the entropy at $p_{00}=p_{01}=0.5$. This
induces a maximum of $\HS_4$ for limit tables which increases with $\lambda$ (see table 1, fourth rows of each odds-ratio). 

\section{Discussion}
In this paper we studied measures of association of 2x2 contingency tables. In contrast to traditional independence analysis, we asked for the 
selection of tables which are far away from independence. This objective was motivated by the analysis of high-dimensional molecular genetic data such
as SNP array data in which a high number of 2x2 tables occur from which one would like to select cases of high dependence called 
{\it linkage disequilibrium}. 

In contrast to detecting a (moderate) deviation from independence, quantifying the strength of association is multiform. 
A large number of possible measures were
proposed in the literature which we shortly reviewed. Many of these measure ($r$, $D^\prime$, $Y$) agree at diagonal tables. Some of the measures become 
extremal for a vanishing diagonal ($r$, $\sMutInf$) while for 
others it suffices that a single cell becomes zero ($D^\prime$, odds-ratio based measures). The measures also markedly differ in cases were one of the 
rows or columns of the table becomes small. Since in practice, it can hardly be decided for these tables whether the dependence is strong or not, 
these tables are not really of interest and are considered as 
{\it junk tables} here. 
Nevertheless, the measure $D^\prime$ can become large in these cases. This is undesirable. $D^\prime$ also varies markedly in a small 
neighbourhood of the vertices of ${\mathbb T}$.    

To study the properties of measures of association, we introduced coordinates on the manifold ${\mathbb T}$ of all tables mapping it to 3-dimensional 
space such that the $x$-axis corresponds to the log-square root of the odds-ratio. We study the measures on the hyperplanes of constant odds-ratio,
looking at the so called {\it margin weighting functions}. These functions are constant for all measures based on the odds-ratio which is 
known to be independent of the margins of the table. For other measures, these functions describe the dependence of 
the measure on the margins for tables with constant odds-ratio. Margin weighting functions illustrate major properties of association measures. 
It helps designing new measures with desired properties, which we demonstrated in the second part of the paper.

The mathematical properties of the margin weighting functions were derived for three measures of association, namely $r$, $\sMutInf$ and $D^\prime$.
It revealed that $r$ and $\sMutInf$ behave very similarly by up-weighting diagonal tables but down-weighting of tables with small rows or columns. In
contrast, $D^\prime$ is not maximal for diagonal tables. Furthermore, it expresses a strange weighting behaviour for tables with small rows and 
columns, up-weighting or down-weighting these tables in dependence on the position of the structural zeros. Such tables occur frequently e.g. 
in SNP data. This property also explains, why the 
estimation problem for $D^\prime$ is not well behaved \citep{SchHas}. 
On the other hand, $D^\prime$ as well as odds-ratio based measures are constructed to up-weight tables which feature a single small entry. These tables represent a prototype of a logical 
table for which one can conclude the state of the column
for one row but not for the other row. These kinds of tables are interesting in genetical statistics since 
they correspond to situations at which no recombinations occurred between two SNPs, i.e. only three of the four theoretically possible haplotypes 
are observed. Therefore, we aimed to define an alternative measure also highlighting L-shaped tables but with a better behaviour at the margins than $D^\prime$ or odds-ratio based 
measures.  

For this purpose, the entropy \citep{Sh} as another canonical structure at 2x2 tables was studied. We proved that the margin weighting function
of this quantity is maximal at the diagonal for odds-ratios within a critical range, namely $\left(W\left(1/e\right)^2,W\left(1/e\right)^{-2}\right)$. 
Outside this range, there are two maxima at L-shaped tables, i.e. tables with a single small cell while the others are (almost) equal. More precisely,
the elements of the opposite diagonal are equal for the maxima. 

The difference between the entropy of a non-diagonal table and the corresponding diagonal table of the same odds-ratio
is a monotone function of the odds-ratio for fixed margins. A new measure of association called $\HS_n$ is defined, which is essentially
{\it Yules Y} weighted by the exponential of this entropy difference. This quantity fulfils all requirements of an association measure, i.e. ranges between -1 and 1, 
is zero in case of independence, is symmetric and a monotone function of the odds-ratio for fixed margins. 
In addition, it agrees with $Y$, $D^\prime$ and $r$ at diagonal tables, up-weights tables with an L-shape and large odds-ratio and is extremal in
case of a single zero in the table. Hence, the measure 
has similar properties than $D^\prime$ except for a better limit behaviour. Since the entropy difference of tables with vanishing row or column is smaller 
than the entropy of the corresponding diagonal table, degenerated tables are markedly down-weighted relative to the diagonal table. The
free constant $n$ allows tuning the degree of this down-weighting. For practical issues we recommend using $n=4$ which 
yields satisfactory results to our experiences. 
However, our procedure of down-weighting junk tables is neither unique nor perfect in the sense that the junk tables are down-weighted 
to zero. The latter one is not possible within the framework of weighting by entropy without loosing other desired properties of the measure,
because the minimum of the absolute differences between the diagonal table and the degenerated tables of the same odds-ratio approaches zero if the 
odds-ratio tends to 0 or $\infty$.  

We recommend using $\HS_4$ instead of $D^\prime$ when interested in selecting L-shaped tables from a large set of tables mostly 
far away from independence and when tables with small marginal frequencies are common. 
When $\HS_4$ is estimated from count data, we recommend using Bayesian plug-in estimators of the frequencies of single cells showing a good compromise 
between accuracy and computational burden \citep{SchHas}.  

\section*{Funding}
This research was funded by the Leipzig Interdisciplinary Research Cluster of Genetic Factors, Clinical Phenotypes, 
and Environment (LIFE Center, University of Leipzig). LIFE is funded by means of the European Union, by the European Regional Development 
Fund (ERDF), the European Social Fund (ESF), and by means of the Free State of Saxony within the framework of its excellence initiative.
{\it Conflict of Interest}: None declared.

\newpage
\section*{Appendix}
\setcounter{equation}{0}
\renewcommand{\theequation}{S.\arabic{equation}}
In this section, we prove the propositions and theorems of our paper. In most situations it is sufficient to consider the case $x>0$ since from 
symmetry conditions the case $x<0$ follows analogously. 
\\
\\
\noindent
{\bf Proposition 6 (Margin weighting function for $r$):} {\it
For all $x \in \mathbb{R}\setminus \left\{0 \right\}$:\\
a) ${r}_{x}$ has exactly one extremum at the origin $\left(y,z\right)=\left(0,0\right)$, corresponding to the diagonal symmetric table with the fixed odds-ratio.\\
b) $\lim_{\left\|(y,z)\right\| \rightarrow \infty }r_x=0$. \\
c) $\lim_{x\to\pm\infty} r_x=\pm 1$ \\
d) $r$ can be extended to $\bar{\mathbb{R}}^3$ except for the lines $(\pm,\pm,*)$ and $(\pm,*,\pm)$ and the vertices $V$. \\
e) $r$ can be extended to $ \overline{{\mathbb T}}$ except for the vertices. }
\\ \\ \noindent
{\it Proof:} 
a) We consider the maximum condition for $r_x$ and $x>0$:
\begin{eqnarray*}
r_x \quad\to\quad\mbox{max.!} &\Leftrightarrow& 
\frac{e^{y+z}}
{\sqrt{\left(e^{x+y+z}+e^y\right)\left(e^{x+y+z}+e^z\right)
\left(e^x+e^y\right)\left(e^x+e^z\right)}} \quad\to\quad\mbox{max.!} \\
&\Leftrightarrow& \left(e^x+e^{-y}\right)\left(e^x+e^{-z}\right)\left(e^x+e^{y}\right)\left(e^x+e^{z}\right) \quad\to\quad\mbox{min.!} \\
&\Leftrightarrow& \left(e^x+e^{-y}\right)\left(e^x+e^{y}\right) \quad\to\quad\mbox{min.!}\quad\wedge\quad y=z \\
&\Leftrightarrow& y=z=0
\end{eqnarray*} 
b) and c) follow easily using equation (3). d) and e) are consequences of b) and c) \hfill $\Box$ 
\\ 
\\ 
\noindent
{\bf Proposition 7 (Margin weighting function for $D^\prime$):}  
{\it For all $x \in \mathbb{R}\setminus \left\{0 \right\}$:\\
a) ${D^\prime}_{x}$ has a non-differentiable edge along the diagonal $y=z$ for $D^\prime>0$ and along the diagonal $y=-z$ for $D^\prime<0$. There is a non-smooth saddle point in the origin.\\
b) 
\begin{eqnarray*}
\lim_{y\to\pm\infty}D^\prime_x&=&\left(e^{2x}-1\right)\cdot\left\{\begin{array}{l@{\quad:\quad}l}
\left(e^{2x}+e^{x\pm z}\right)^{-1} & x>0 \\
\left(e^{x\mp z}+1\right)^{-1} & x<0 
\end{array}\right. \\ 
\lim_{z\to\pm\infty}D^\prime_x&=&\left(e^{2x}-1\right)\cdot\left\{\begin{array}{l@{\quad:\quad}l}
\left(e^{2x}+e^{x\pm y}\right)^{-1} & x>0 \\
\left(e^{x\mp y}+1\right)^{-1} & x<0 
\end{array}\right. 
\end{eqnarray*}
Thus, limit functions have a range of $\left(0,1-e^{-2x}\right)$ for $x>0$ and 
$\left(e^{2x}-1,0\right)$ for $x<0$, where 0 is obtained for 
$y\to\pm\infty$, $z\to\pm\infty$, $x>0$ and $y\to\mp\infty$, $z\to\pm\infty$, $x<0$.\\
c) $\lim_{x\to\pm\infty} D^\prime_x=\pm 1$ \\
d) $D^\prime$ can be extended to $\bar{\mathbb{R}}^3$ except for the vertices $V_g$. \\
e) $D^\prime$ can be extended to $ \overline{{\mathbb T}}$ except for the edges and vertices. }
\\ \\ \noindent
{\it Proof:} 
a) Assume $x>0$, consider the path $y=w+c$, $z=w-c$, $w=\;$const. Calculating the left-hand and right-hand derivative of $D_{max}$ at $c=0$ using equation (4) yields:
\begin{eqnarray*}
\lim_{c\to 0\pm 0}\frac{d}{dc} D_{max}&=&\mp e^w\left(e^x+2e^w+e^{x+2w}\right)
\end{eqnarray*}
Since the term in parentheses is positive, $D_{max}$ has a wedge at $c=0$. On the other hand, $D^\prime_x$ has a maximum at $y=z=0$ along the path $y=z$ since
\begin{eqnarray*}
D^\prime_x \left(y,y\right) \quad\to\quad\mbox{max.!} &\Leftrightarrow& 
\frac{e^{2y}}{\left(e^{x+2y}+e^y\right)\left(e^x+e^y\right)} \quad\to\quad\mbox{max.!} \\
 &\Leftrightarrow& \left(e^x+e^{-y}\right)\left(e^x+e^{y}\right) \quad\to\quad\mbox{min.!} \\
&\Leftrightarrow& y=0
\end{eqnarray*}
Hence, $D^\prime$ has a non-differentiable saddle point at $y=z=0$. The case $x<0$ follows analogously. The limit behaviour considered in b) to e) is easy to see using equation (4). \hfill $\Box$
\\ 
\\ \noindent
{\bf Proposition 8 (Margin weighting function for $\sMutInf$):}   
{\it For all $x \in \mathbb{R}\setminus \left\{0 \right\}$:\\
a) $\sMutInf_{x}$ has exactly one maximum at the origin $\left(y,z\right)=\left(0,0\right)$.\\
b) $\lim_{\left\|(y,z)\right\| \rightarrow \infty }\sMutInf_x=0$. \\
c) $\lim_{x\to\pm\infty} \sMutInf_x=\pm\left(\log_2\left( e^{y\pm z}+1\right)-\frac{e^{y\pm z}}{e^{y\pm z}+1}\log_2 e^{y \pm z}  \right)$\\ 
Thus $\sMutInf_{x}\to \pm 1$ for $y=\mp z$ and $x \to\pm\infty$ respectively. \\
d) $\sMutInf$ can be extended to $\bar{\mathbb{R}}^3$ except for the vertices $V_b$. \\
e) $\sMutInf$ can be extended completely to $ \overline{{\mathbb T}}$. }
\\ \\ \noindent
{\it Proof:} 
a) We consider tables 
$t_\mu=\frac{1}{N_\mu}\left({\mu p_{00} \atop p_{10}}{p_{01} \atop p_{11}/\mu}\right)$ and 
$t_\nu=\frac{1}{N_\nu}\left({ p_{00} \atop \nu p_{10}}{p_{01}/\nu \atop p_{11}}\right)$
 of the same odds-ratio than $t$ for $\mu,\nu>0 $ and $N_\mu$ and $N_\nu$ are the normalisation constants
$N_\mu=\mu p_{00} + p_{01} + p_{10} + p_{11}/\mu$ and $N_\nu= p_{00} + \nu p_{01} + p_{10}/\nu + p_{11}$ respectively.
We aim to proof that 
\begin{eqnarray}
\left.\frac{d}{d\mu}{\sMutInf}_x\left(t_\mu\right)\right|_{\mu=1}=0 &\Leftrightarrow& p_{00}=p_{11} \\
\left.\frac{d}{d\nu}{\sMutInf}_x\left(t_\nu\right)\right|_{\nu=1}=0 &\Leftrightarrow& p_{01}=p_{10}
\end{eqnarray}
Assuming $\lambda>1$, $p_{00}\le p_{11}$ without restriction of generality, we obtain after some calculations
\begin{eqnarray}
\left.\frac{d}{d\mu}{\sMutInf}_x\left(t_\mu\right)\right|_{\mu=1}&=&
\left(p_{11}-p_{00}\right){\sMutInf}_x\left(t\right) \\ \nonumber
&&{}+p_{00}\log_2\frac{p_{00}}{\left(p_{00}+p_{01}\right)\left(p_{00}+p_{10}\right)} -p_{11}\log_2
\frac{p_{11}}{\left(p_{11}+p_{01}\right)\left(p_{11}+p_{10}\right)} \\ \nonumber
&=&\left(p_{11}-p_{00}\right){\sMutInf}_x\left(t\right) \\ \nonumber
&&{}+p_{00}p_{11}\left(\frac{\log_2\left(1+p_{00}\left(\frac{1}{\lambda}-1\right)\right)}{p_{00}}-
\frac{\log_2\left(1+p_{11}\left(\frac{1}{\lambda}-1\right)\right)}{p_{11}}\right)
\end{eqnarray}
The first term is non-negative and equals zero iff $p_{00}=p_{11}$. 
\\ \noindent
Consider the monotonicity of the term $\frac{\log_2\left(1+x\left(\frac{1}{\lambda}-1\right)\right)}{x}$ for
$x\in\left(0,1\right)$. It holds that
\begin{eqnarray}
\frac{d}{dx}\frac{\log_2\left(1+x\left(\frac{1}{\lambda}-1\right)\right)}{x} &=&\frac{1}{x^2\ln2}\left(
-\ln z+\frac{z-1}{z}\right)
\end{eqnarray}
with $z=x\left(\frac{1}{\lambda}-1\right)+1\in \left(0,1\right)$. In this interval, (A.4) is negative since
$-\ln z+\frac{z-1}{z}$ is monotonically increasing for $z\in\left(0,1\right)$, taking its maximum for $z\to 1$.
Hence $\frac{\log_2\left(1+x\left(\frac{1}{\lambda}-1\right)\right)}{x}$ is monotonically decreasing for 
$x\in\left(0,1\right)$. In conclusion, the second term of (A.3) is non-negative too and equals zero iff $p_{00}=p_{11}$.
This proves (A.1). Analogously, using $t_\nu$ instead of $t_\mu$ proves (A.2). 
\\ \noindent
b)-e) are obvious exploiting the continuity of the functions involved, i.e. taking the limit of the tables first. \hfill $\Box$
\\ 
\\ 
\noindent 
{\bf Theorem 1 (magic odds-ratio):}
{\it Define the "magic odds-ratio" by $L_{magic}=W(1/e)^{-2} \approx 12.89 $. 
Let $L > 1$. The entropy $H$ restricted to the submanifold of constant odds-ratio $L$ in $\mathbb{T}$
\begin{itemize}
	\item has a single maximum at the diagonal table of odds-ratio $L$ if $1 < L \le L_{magic}$.

	\item has a saddle point at the diagonal table of odds-ratio $L$ and two "L-shaped" tables as maxima which transpose with matrix transposition if $L_{magic} < L$.
\end{itemize}
"L-shaped" means that for $L \rightarrow \infty$ one of the maxima approaches the table $\left( {1/3\atop 0}{1/3\atop 1/3}\right)$.  
For the case $L<1$ a similar result can be derived by transposing principal and secondary diagonals. }
\\ \\ \noindent
{\it Proof:}
The constraint odds-ratio = $L$ can be written in the form:
\begin{eqnarray}
 \ln(p_{00}) - \ln(p_{01}) - \ln(p_{10}) +\ln(p_{11})= \ln(L) 
\end{eqnarray}
We assume $L>1$ in the following without restriction of generality since the case $L<1$ can be studied analogously. 
A second constraint is given by
\begin{eqnarray}
 p_{00} + p_{01} + p_{10} + p_{11} = 1 
\end{eqnarray}
In order to study the critical points of $H$, we now consider the extremal value problem of $H$ given the constraints (A.5) and (A.6). For this
purpose, we introduce Lagrange multipliers $\Lambda_1$ and $\Lambda_2$ and determine the first variation of the following function:
\begin{eqnarray}
f\left(t,\Lambda_1,\Lambda_2\right) \gl  
-\left(p_{00} \cdot \ln(p_{00}) +p_{01} \cdot \ln(p_{01}) +p_{10} \cdot \ln(p_{10}) + p_{11} \cdot \ln(p_{11})\right) + \\ \nonumber
{}+\Lambda_1 \cdot  \left( \ln(p_{00}) - \ln(p_{01}) - \ln(p_{10}) +\ln(p_{11}) -\ln(L)\right) + \Lambda_2 \cdot \left( p_{00} +p_{01}+p_{10}+p_{11}-1 \right) 
\end{eqnarray}
Calculating the partial derivatives $\frac{\partial}{\partial p_{ij}}$ gives four equations:
\begin{eqnarray*}
\ln(p_{00}) +1 - \frac{\Lambda_1}{p_{00}} -\Lambda_2 &=&0\\
\ln(p_{11}) +1 - \frac{\Lambda_1}{p_{11}} -\Lambda_2 &=&0\\
\ln(p_{10}) +1 + \frac{\Lambda_1}{p_{10}} -\Lambda_2 &=&0\\
\ln(p_{01}) +1 + \frac{\Lambda_1}{p_{01}} -\Lambda_2 &=&0
\end{eqnarray*}
In order to solve this system explicitly, we recall that Lambert's W function is defined\\ as the inverse function to $x \exp(x)$. Hence, it holds that
\begin{eqnarray*}
p_{ij}=\frac{\pm \Lambda_1}{W(\pm \Lambda_1 \exp(1- \Lambda_2 ))}
\end{eqnarray*}
where the upper sign holds for $p_{00}$ and $p_{11}$ and the lower sign for $p_{01}$ and $p_{10}$ respectively. At the first look it seems as that the only solution is the diagonal-symmetric table. 
But $W$ is a multi-branch function since $y=x \exp(x)$ has two solutions for $y \in (-1/e,0)$. The two real-valued braches are traditionally called 
$W_0$ when $x \in (-1/e,0)$ and $W_{-1}$ when $x \in (-\infty,-1/e)$. Note that $W_{-1}\left(-1/e\right)=W_0\left(-1/e\right)=-1$.\\    
\noindent
Assume $\Lambda_1>0$. Inserting the solutions for $p_{ij}$ in the condition on $\ln(L)$ we get three possible solutions.\\
\noindent
a) $\ln(L)=\ln(\frac{W_0(-z)^2}{W(z)^2})$ \\
This solution exists only for $L \in \left(1,W(1/e)^{-2}\right]$. $W(1/e)^{-2} \approx 12.89615...$. \\
\noindent
b) $\ln(L)=\ln(\frac{W_{-1}(-z)^2}{W(z)^2})$ \\
This solution exists only for $L \in \left[W(1/e)^{-2},\infty\right)$. \\
\noindent
c) $\ln(L)=\ln(\frac{W_{-1}(-z)W_0(-z)}{W(z)^2})$\\ 
These solutions exist only for $L \in \left[W(1/e)^{-2},\infty\right)$. \\
\noindent
Hence, we have a single critical point for $L \in \left(1,W(1/e)^{-2}\right)$ but three critical points for $L \in \left(W(1/e)^{-2},\infty\right)$. The next lemma characterises these critical points.
\\ \\ \noindent
{\bf Lemma: (Characterisation of the critical points of $H$ for given odds-ratio):}\\ \noindent {\it 
a) For $L \in \left(1,W(1/e)^{-2}\right]$, $H$ has a maximum at the diagonal table of odds-ratio $L$. \\
b) For $L \in \left(W(1/e)^{-2},\infty\right)$, $H$ has a saddle-point at the diagonal table and two maxima at the other two critical points. If $L\to\infty$ these maxima tend to the tables
$\left({1/3 \atop 0 }{1/3 \atop 1/3} \right)$ and $\left({1/3 \atop 1/3 }{0 \atop 1/3} \right)$ respectively. }
\\ \\ \noindent
{\it Proof:} We study the following tables: 
$t_\mu=\frac{1}{N_\mu}\left({\mu a \atop c}{c \atop a/\mu}\right)$ with $N_\mu=\mu a + \frac{a}{\mu} +2c$ and $t_\nu=\frac{1}{N_\nu}\left({a \atop c/\nu}{\nu c \atop a}\right)$ 
with $N_\nu=2a + \mu c + \frac{c}{\nu}$, $a,c>0$, $a+c=1/2$. We calculate the second derivative of $H\left(t_\mu\right)$ and $H\left(t_\nu\right)$ at $\mu=1$ and $\nu=1$ respectively.  
After some calculations one obtains
\begin{eqnarray*}
\left.\frac{d^2H}{d\mu^2}\left(t_\mu\right)\right|_{\mu=1}&=&-\frac{1}{\ln 2}\frac{\sqrt{\lambda}}{1+\sqrt{\lambda}}\left(1+\frac{1}{1+\sqrt{\lambda}}\ln\sqrt{\lambda}\right)<0
\end{eqnarray*}
In contrast
\begin{eqnarray*}
\left.\frac{d^2H}{d\nu^2}\left(t_\nu\right)\right|_{\nu=1}&=&-\frac{1}{\ln 2}\frac{1}{1+\sqrt{\lambda}}\left(1-\frac{\sqrt{\lambda}}{1+\sqrt{\lambda}}\ln\sqrt{\lambda}\right)
\end{eqnarray*}
Which is greater than 0 for $\lambda>W\left(1/e\right)^{-2}$. Thus $\left({a \atop c}{c \atop a}\right)$ becames a saddle point for $\lambda>W\left(1/e\right)^{-2}$ but is a maximum for
$\lambda \le W\left(1/e\right)^{-2}$. The other suppositions of the lemma and theorem 1 are then easy to see. \hfill$\Box$
\\ 
\\
\noindent
{\bf Lemma 1 (Monotony of the entropy difference):} {\it Let $H$ be the entropy of $t$ and $H_{diag}$ be the entropy of the corresponding diagonal table of the same odds-ratio $\lambda$. 
Then, $H_{diag}-H$ is monotonically decreasing for increasing $\lambda>1$ and constant margins. }
\\ \\ \noindent
{\it Proof:}
Let $\ve>0$ and $t_\ve=\left({p_{00}+\ve\atop p_{10}-\ve}{p_{01}-\ve\atop p_{11}+\ve}\right)$ a table 
with increased odds-ratio but same margins compared to $t$. We show that
$\left.\frac{d}{d\ve}\right|_{\ve=0}H_{diag}\left(t_\ve\right)-H\left(t_\ve\right)\le 0$, where equality holds iff $t$ is diagonal. 
After some calculations we obtain
\begin{eqnarray}
\left.\frac{d}{d\ve}\right|_{\ve=0}H_{diag}\left(t_\ve\right) &=&-\frac{\sqrt{\lambda}\log_2\lambda}{4\left(1+\sqrt{\lambda}\right)^2}\sum_{i,j=0}^1
\frac{1}{p_{ij}} \\
\left.\frac{d}{d\ve}\right|_{\ve=0}H\left(t_\ve\right) &=& -\log_2\lambda
\end{eqnarray}
Thus
\begin{eqnarray}
\left.\frac{d}{d\ve}\right|_{\ve=0}\left(H_{diag}\left(t_\ve\right)- H\left(t_\ve\right) \right) &=&
\frac{\log_2\lambda}{4}\left(4-\frac{\sqrt{\lambda}}{\left(1+\sqrt{\lambda}\right)^2}\sum_{i,j=0}^1\frac{1}{p_{ij}}\right)
\end{eqnarray}
Now consider the tables $t_\mu=\frac{1}{N_\mu}\left({\mu p_{00} \atop p_{10}}{p_{01} \atop p_{11}/\mu}\right)$ and 
$t_\nu=\frac{1}{N_\nu}\left({ p_{00} \atop \nu p_{10}}{p_{01}/\nu \atop p_{11}}\right)$ 
of the same odds-ratio than $t$ for $\mu,\nu\ge 1$ and the normalisation constants $N_\mu=\mu p_{00} + p_{01} + p_{10} + p_{11}/\mu$ and
$N_\nu= p_{00} + \nu p_{01} + p_{10}/\nu + p_{11}$ respectively.
\\
Assume $p_{00}\le p_{11}$ and $p_{01}\le p_{10}$ without restriction of generality, we see that for $f\left(t\right)=\sum_{i,j=0}^1\frac{1}{p_{ij}}$
it holds that
\begin{eqnarray*}
\left.\frac{d}{d\mu}\right|_{\mu=1} f\left(t_\mu\right) &=&\left(p_{00}-p_{11}\right)\left(\sum_{i,j=0}^1\frac{1}{p_{ij}}+\frac{1}{p_{00}p_{11}}\right)\le 0 \\
\left.\frac{d}{d\nu}\right|_{\nu=1} f\left(t_\nu\right) &=&\left(p_{10}-p_{01}\right)\left(\sum_{i,j=0}^1\frac{1}{p_{ij}}+\frac{1}{p_{01}p_{10}}\right)\le 0 
\end{eqnarray*}
were equality holds iff $t$ is diagonal. Hence the maximum of the term in parenthesis of (A.10) is obtained iff 
$t$ is diagonal. On the other hand, for $t$ diagonal it hold that
\begin{eqnarray}
\frac{\sqrt{\lambda}}{\left(1+\sqrt{\lambda}\right)^2}\sum_{i,j=0}^1\frac{1}{p_{ij}}-4&=&0
\end{eqnarray}
\hfill $\Box$

\newpage
\section*{Tables and Figures}
\begin{tabular}{||l|l|l|l||l||l|l|l|l||} \hline
$p_{00}$ & $p_{01}$ & $p_{10}$ & $p_{11}$ & $\lambda$ & $Y$ & $r$ & $D^\prime$ & $\HS_4$ \\ \hline\hline
0.25    & 0.25   & 0.25    & 0.25    & 1 & 0     & 0     & 0       & 0 \\ \hline
0.25    & 0.25   & 0.25    & 0.25    & 1 & 0     & 0     & 0       & 0 \\ \hline
{\it 0} & 1      & {\it 0} & {\it 0} & 1 & 0     & 0     & 0       & 0 \\ \hline
0.5     & 0.5    & {\it 0} & {\it 0} & 1 & 0     & 0     & 0       & 0 \\ \hline
1       & {\it 0}& {\it 0} & {\it 0} & 1 & 0     & 0     & 0       & 0 \\ \hline\hline
0.293   & 0.207  & 0.207   & 0.293   & 2 & 0.172 & 0.172 & 0.172   &0.172 \\ \hline
0.286   & 0.286  & 0.143   & 0.286   & 2 & 0.172 & 0.167 & 0.222   &0.139 \\ \hline
{\it 0} & 1      & {\it 0} & {\it 0} & 2 & 0.172 &{\it 0}& 0.5     &{\it 0} \\ \hline
0.5     & 0.5    & {\it 0} & {\it 0} & 2 & 0.172 &0.002  & 0.333   &{\it 0} \\ \hline
1       &{\it 0} & {\it 0} & {\it 0} & 2 & 0.172 &{\it 0}& {\it 0} &{\it 0} \\ \hline\hline
0.345   &0.155   &0.155    &0.345    & 5 &0.382  &0.382  &0.382    &0.382 \\ \hline
0.312   &0.312   &0.062    &0.312    & 5 &0.382  &0.333  &0.556    &0.282 \\ \hline
{\it 0} &1       &{\it 0}  &{\it 0}  & 5 &0.382  &{\it 0}&0.8      &{\it 0} \\ \hline
0.5     &0.5     &{\it 0}  &{\it 0}  & 5 &0.382  &0.005  &0.667    &{\it 0}\\ \hline
1       &{\it 0} &{\it 0}  &{\it 0}  & 5 &0.382  &{\it 0}&{\it 0}  &{\it 0}\\ \hline\hline
0.38    &0.12    &0.12     &0.38     &10 &0.519  &0.519  &0.519    &0.519\\ \hline
0.323   &0.323   &0.032    &0.323    &10 &0.519  &0.409  &0.744    &0.441\\ \hline
{\it 0} &1       &{\it 0}  &{\it 0}  &10 &0.519  &{\it 0}&0.9      &{\it 0}\\ \hline
0.5     &0.5     &{\it 0}  &{\it 0}  &10 &0.519  &0.007  &0.818    &{\it 0}\\ \hline
1       &{\it 0} &{\it 0}  &{\it 0}  &10 &0.519  &{\it 0}&{\it 0}  &{\it 0}\\ \hline\hline
0.409   &0.091   &0.091    &0.409    &20 &0.635  &0.635  &0.635    &0.635\\ \hline
0.328   &0.328   &0.016    &0.328    &20 &0.635  &0.452  &0.862    &0.627\\ \hline
{\it 0} &1       &{\it 0}  &{\it 0}  &20 &0.635  &{\it 0}&0.95     &{\it 0}\\ \hline
0.5     &0.5     &{\it 0}  &{\it 0}  &20 &0.635  &0.009  &0.905    &0.001\\ \hline
1       &{\it 0} &{\it 0}  &{\it 0}  &20 &0.635  &{\it 0}&{\it 0}  &{\it 0}\\ \hline\hline
0.438   &0.062   &0.062    &0.438    &50 &0.752  &0.752  &0.752    &0.752\\ \hline
0.331   &0.331   &0.007    &0.331    &50 &0.752  &0.48   &0.942    &0.821\\ \hline
{\it 0} &0.999   &{\it 0}  &{\it 0}  &50 &0.752  &{\it 0}&0.98     &{\it 0}\\ \hline
0.5     &0.5     &{\it 0}  &{\it 0}  &50 &0.752  &0.012  &0.961    &0.086\\ \hline
1       &{\it 0} &{\it 0}  &{\it 0}  &50 &0.752  &{\it 0}&{\it 0}  &{\it 0}\\ \hline\hline
0.455   &0.045   &0.045    &0.455    &100&0.818  &0.818  &0.818    &0.818\\ \hline
0.332   &0.332   &0.003    &0.332    &100&0.818  &0.49   &0.97     &0.904\\ \hline
{\it 0} &0.999   &{\it 0}  &{\it 0}  &100&0.818  &{\it 0}&0.99     &{\it 0}\\ \hline
0.5     &0.5     &{\it 0}  &{\it 0}  &100&0.818  &0.015  &0.98     &0.316\\ \hline
1       &{\it 0} &{\it 0}  &{\it 0}  &100&0.818  &{\it 0}&{\it 0}  &{\it 0}\\ \hline\hline
\end{tabular}
\\ \noindent
{\bf Table 1 (Measures of association for selected tables):} Tables entries rounded to three decimals are presented in 
columns 1 to 4. Normally printed zeros are hard zeros
while zeros in italic are values less than 0.0005
The fifth columns presents the odds-ratio of the tables. For each odds-ratio, we studied five tables: the diagonal table (first row of the corresponding odds-ratio), 
a table with three equal entries (second row), a table for which it holds that $p_{01}\approx 1$ or $p_{00}\approx 1$ (third and fifth row respectively) and a table for 
which $p_{00}\approx p_{01}\approx 0.5$ (fourth row). The last four columns contain the corresponding measures of association rounded to three decimals. 
.

\newpage
\scalebox{1}{\includegraphics{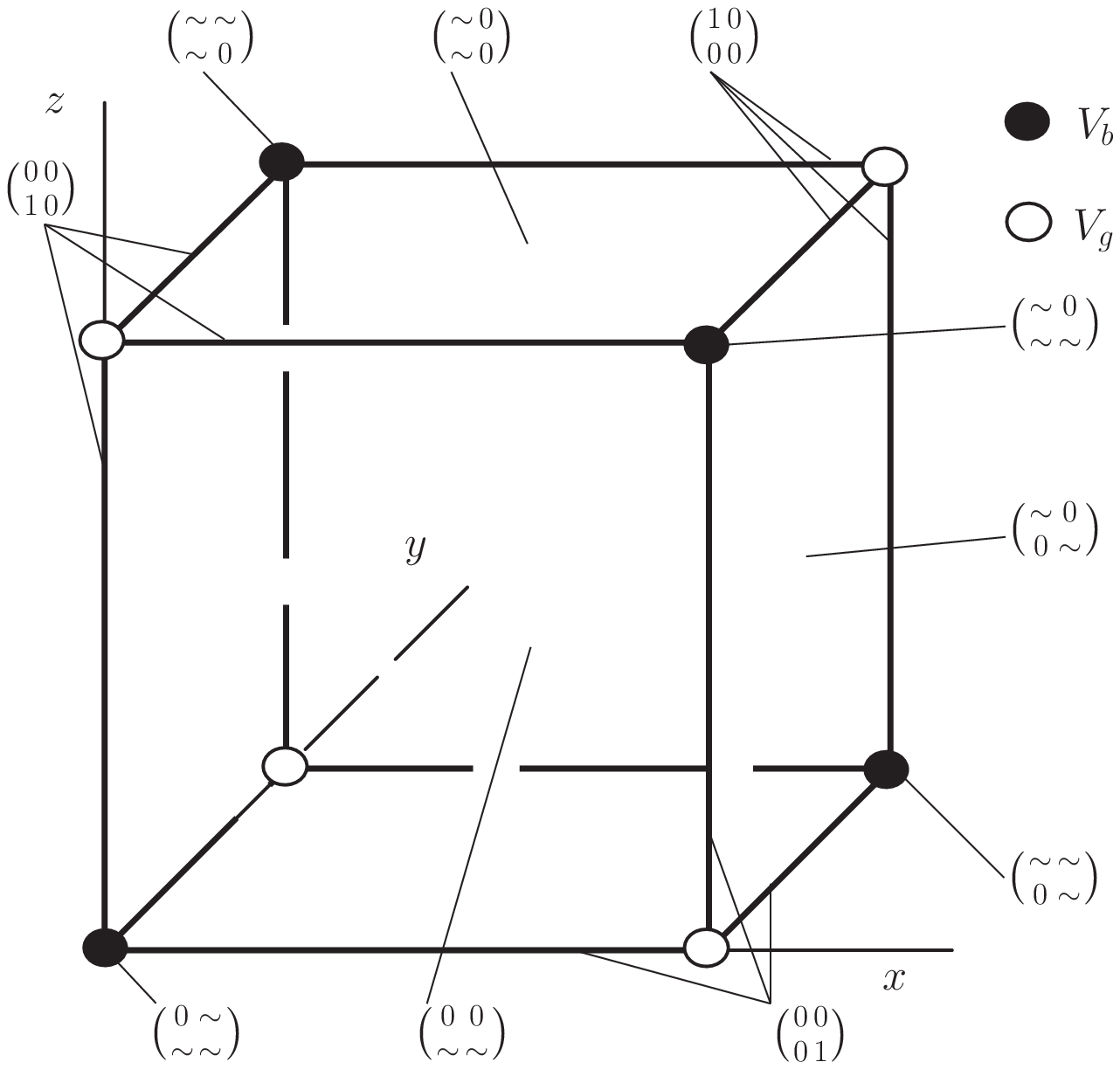}} 
\\ \noindent
{\bf Figure 1:} Illustration of the maps $\Theta$ and $\Psi$ on the boundaries of $\bar{\mathbb{R}}^3$ and $\overline{{\mathbb T}}$. ''$\sim$'' represents 
positive number adding up to 1.

\newpage
\noindent
\scalebox{0.67}{\includegraphics{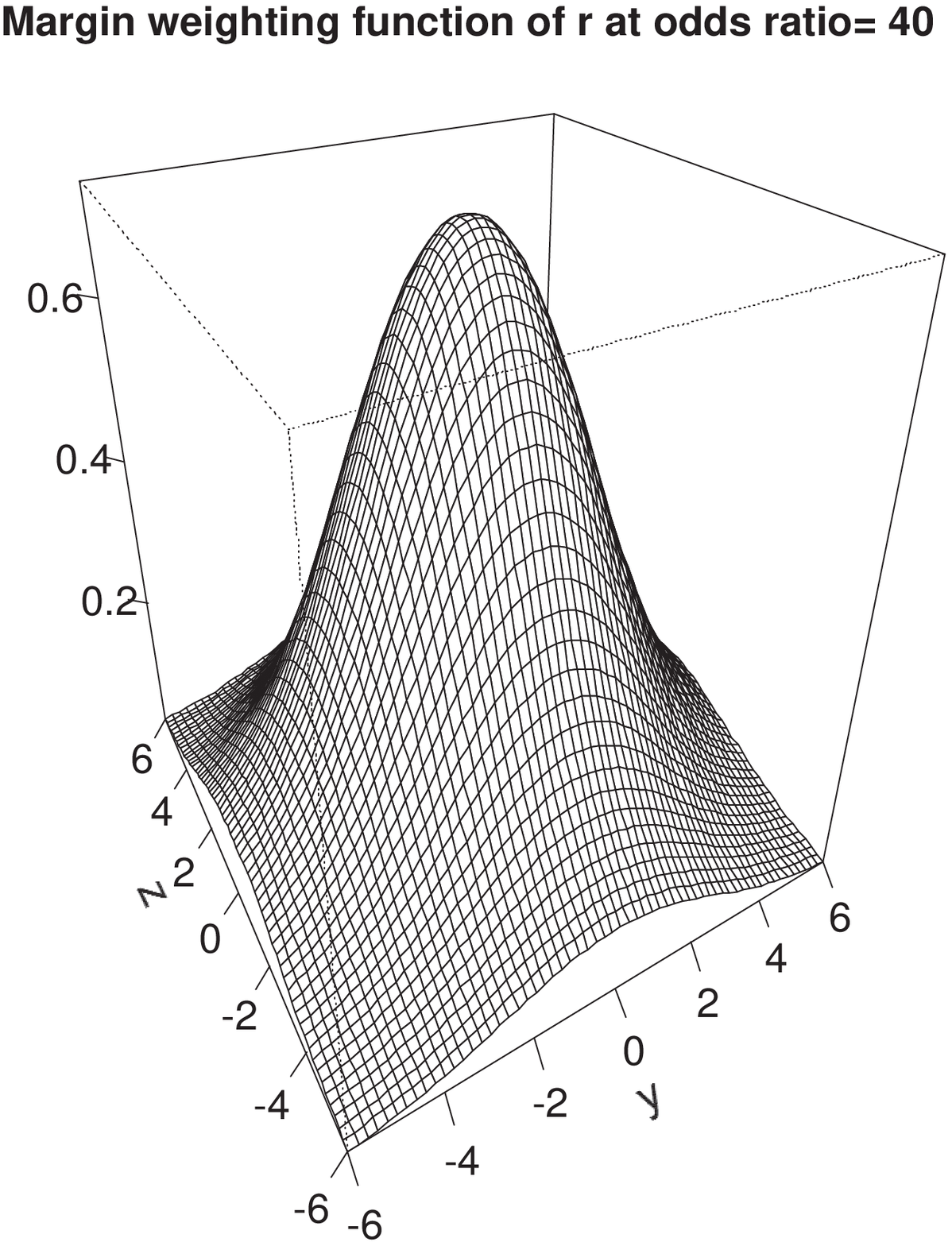}}
\\ \noindent
{\bf Figure 2:} Margin-weighting function of $r$

\newpage
\noindent
\scalebox{0.67}{\includegraphics{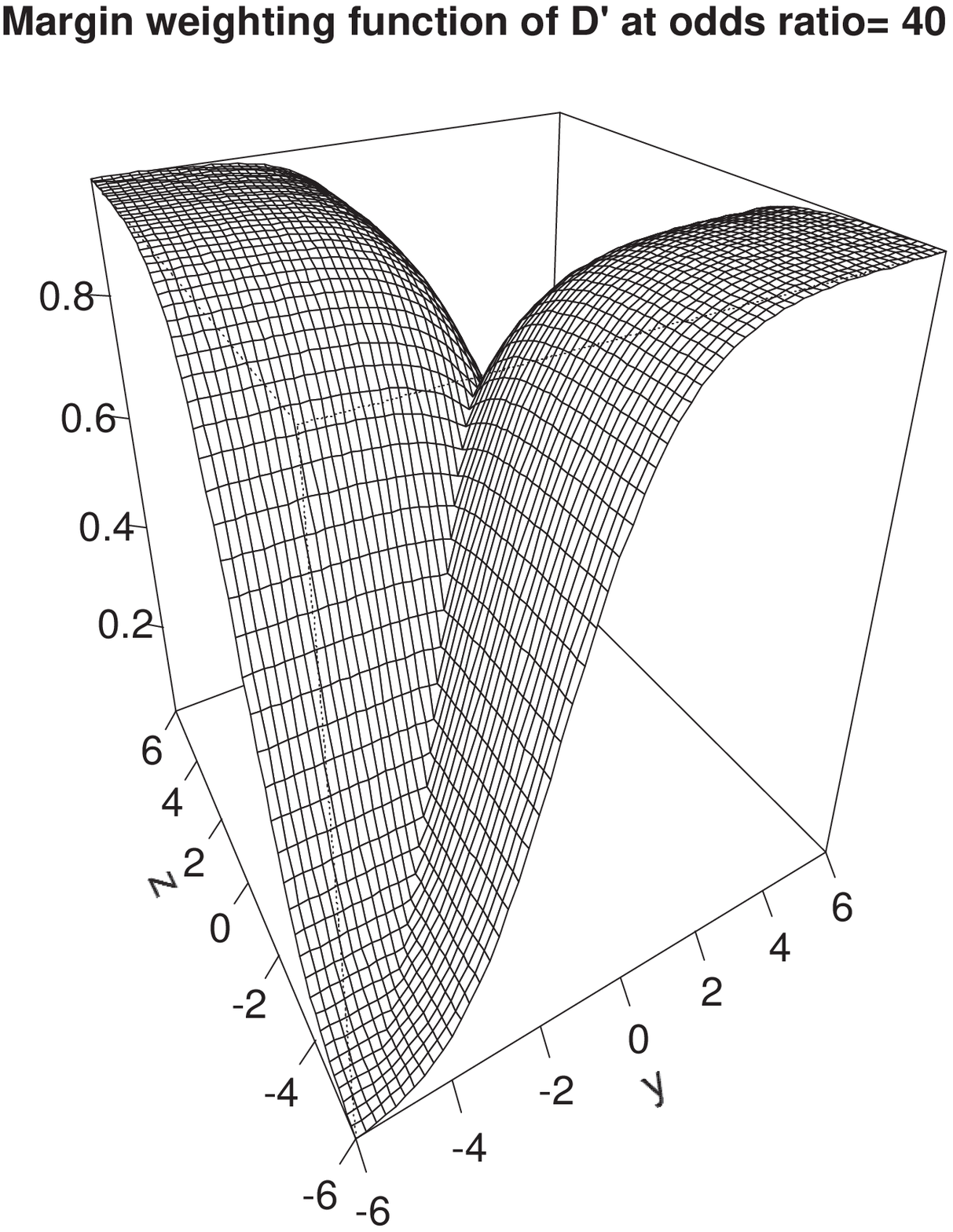}}
\\ \noindent
{\bf Figure 3:} Margin-weighting function of $D^\prime$

\newpage
\noindent
\scalebox{0.67}{\includegraphics{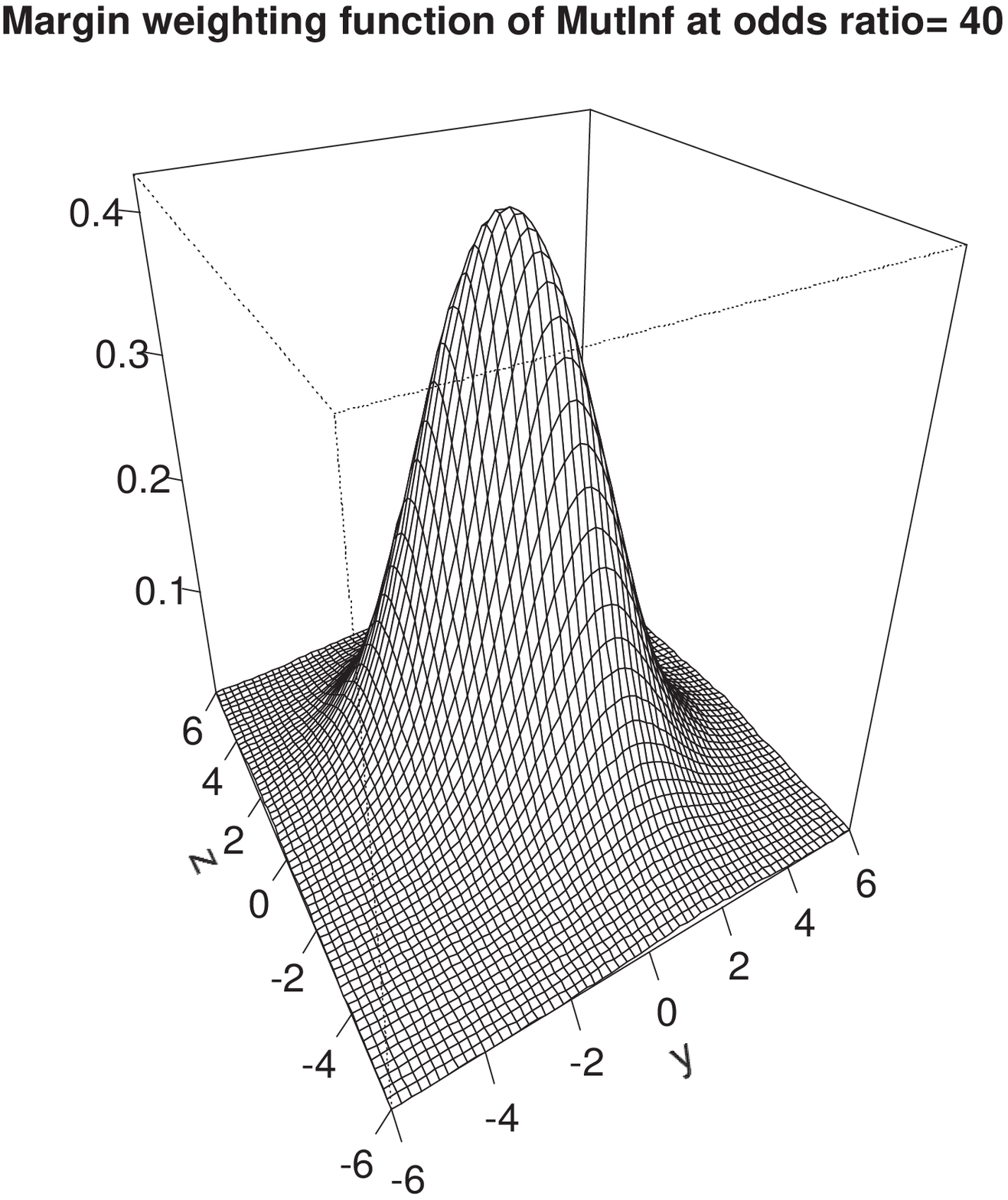}}
\\ \noindent
{\bf Figure 4:} Margin-weighting function of $\sMutInf$

\newpage
\noindent
\scalebox{0.5}{\includegraphics{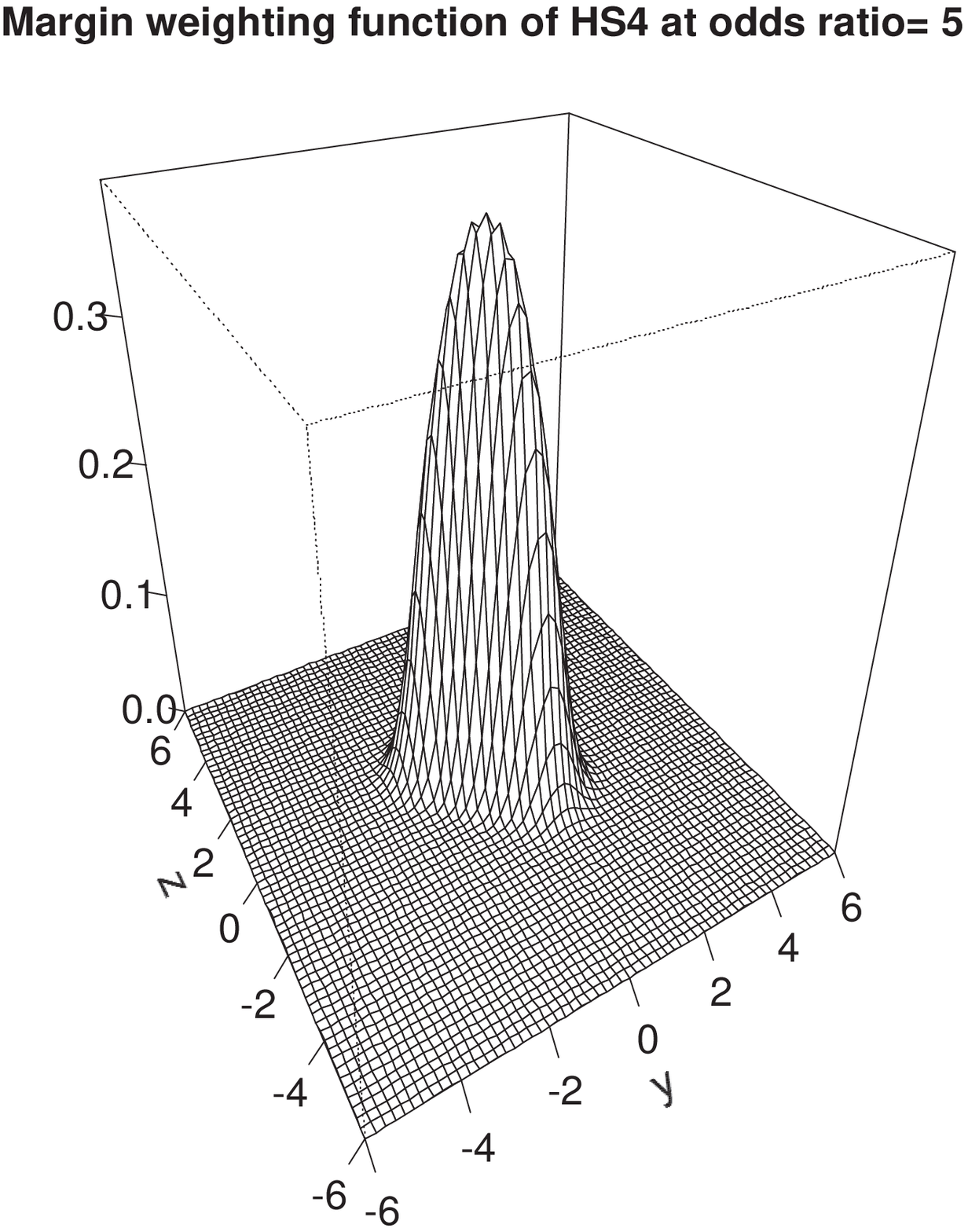}}
\scalebox{0.5}{\includegraphics{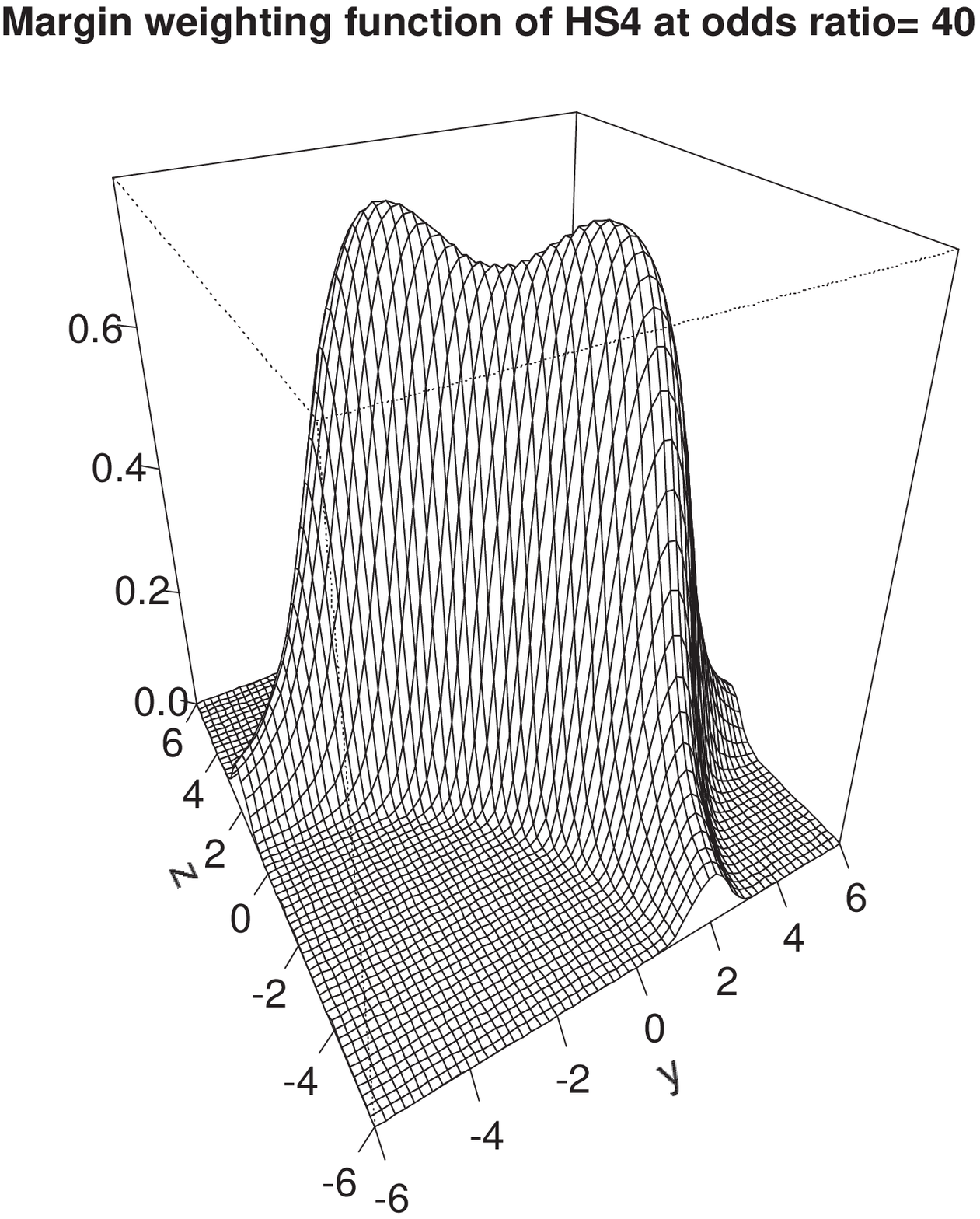}} \\ \noindent
{\bf Figure 5:} Margin-weighting function of $\HS_4$

\end{document}